\documentclass[12pt]{article}
\usepackage{framed}
\RequirePackage{kpfonts}
\RequirePackage[sf,bf,small,raggedright]{titlesec}
\usepackage{textcomp}
\usepackage{amssymb}
\usepackage{bbm}
\usepackage{fancyhdr}
\usepackage{amsmath}
\usepackage{amsbsy,amsthm}
\usepackage{amscd}
\usepackage{latexsym}
\usepackage{graphicx}   
\usepackage{pdfsync}
\usepackage{blkarray}
\usepackage{multirow}
\usepackage{hyperref}
\usepackage{tcolorbox}

\newcommand{\R}{\mathbb{R}}

\newcommand{\inr}[1]{\left\langle #1 \right\rangle}

\newcommand{\E}{\mathbb{E}}

\newcommand{\eps}{\varepsilon}

\newtheorem{theorem}{Theorem}
\newtheorem{corollary}{Corollary}

\newtheorem{question}{Question}
\newtheorem{remark}{remark}

\newtheorem{example}{Example}

\numberwithin{equation}{section}

\def \proof {\noindent {\bf Proof.}\ \ }

\def \endproof
{{\mbox{}\nolinebreak\hfill\rule{2mm}{2mm}\par\medbreak}}
\def\IND{\mathbbm{1}}

\newcommand{\ol}{\overline}
\newcommand{\wt}{\widetilde}
\newcommand{\wh}{\widehat}

\newcommand{\EXP}{\mathbb{E}}
\newcommand{\PROB}{\mathbb{P}}

\newcommand{\Tr}{\mathrm{Tr}}

\newcommand{\defeq}{\stackrel{\mathrm{def.}}{=}}

\begin{document}

\title{Near-optimal mean estimators with respect to general norms
\thanks{
G\'abor Lugosi was supported by
the Spanish Ministry of Economy and Competitiveness,
Grant MTM2015-67304-P and FEDER, EU. Shahar Mendelson was supported in part by the Israel Science Foundation.
}
}
\author{
G\'abor Lugosi\thanks{Department of Economics and Business, Pompeu
  Fabra University, Barcelona, Spain, gabor.lugosi@upf.edu}
\thanks{ICREA, Pg. Lluís Companys 23, 08010 Barcelona, Spain}
\thanks{Barcelona Graduate School of Economics}
\and
Shahar Mendelson \thanks{Mathematical Sciences Institute, The Australian National University and Department of Mathematics, Technion, I.I.T,  shahar.mendelson@anu.edu.au}}

\maketitle

\begin{abstract}
  We study the problem of estimating the mean of a random vector in
  $\R^d$ based on an i.i.d.\ sample, when the accuracy of the
  estimator is measured by a general norm on $\R^d$. We construct an
  estimator (that depends on the norm) that achieves an essentially
  optimal accuracy/confidence tradeoff under the only assumption that
  the random vector has a well-defined covariance matrix. The
  estimator is based on the construction of a uniform median-of-means
  estimator in a class of real valued functions that may be of
  independent interest.
\end{abstract}

\section{Introduction}
In this note we explore the problem of multivariate mean estimation
with respect to an arbitrary norm. To formulate the question, let $\|
\cdot \|$ be a norm on $\R^d$ and let $X$ be a random vector in
$\R^d$. One only assumes that $X$ has a mean $\mu = \E X$ and a
well-defined covariance matrix $\Sigma=\EXP (X-\mu)\otimes(X-\mu)$.
The statistical problem we consider is estimating the mean vector $\mu$ from a sample $(X_i)_{i=1}^N$ of $N$ independent copies of $X$.
We do not assume any knowledge on the distribution. The goal is to approximate the mean $\mu$ by finding
some \emph{mean estimator}
$\wh{\mu}_N = \wh{\mu}_N(X_1,\ldots,X_N)\in \R^d$ such that $\|\wh{\mu}_N-\mu\|$ is as small as possible.

Formally, the problem studied in this note is as follows:

\begin{tcolorbox}
Given a norm $\| \cdot \|$, a confidence parameter $\delta\in (0,1)$ and an i.i.d.\ sample of cardinality $N$, find an estimator
$\wh{\mu}_N$ and the best possible accuracy $\eps$ for which
$$
\|\wh{\mu}_N - \mu \| \leq \epsilon \ \ \ {\rm with \ probability \ at \ least \ } 1-\delta~.
$$
\end{tcolorbox}

Various versions of this question have been studied extensively in recent years, but it was far from resolved. In fact, even the correct order of the best accuracy $\epsilon$ was not clear, except in special situations.
While there are some results for specific choices of norms, the only estimate that is known to be optimal was obtained in Lugosi and Mendelson \cite{LuMe16a} for the Euclidean norm, see also Joly, Lugosi, and Oliveira \cite{JoLuOl17} and Catoni and Giulini \cite{CaGi17}.
In addition, there are also several partial results (see Minsker \cite{Min16}, Catoni and Giulini \cite{CaGi17}) for other special norms (mainly in the context of the matrix operator norm) and which are suboptimal, as we will see below.

We start by discussing what kind of accuracy $\epsilon$ one should be aiming for. To this end, first consider the case when $X$ is a real-valued random variable with finite mean $\mu$ and variance $\sigma^2$. Since the real-valued case is well-understood, it will eventually lead
us to the possible identity of $\epsilon$ in the vector-valued scenario.

The first observation (see, e.g., Catoni \cite{Cat10}) is that if $X$ is a Gaussian random variable then the best mean estimate
that one can hope for is such that, with probability $1-\delta$,
\begin{equation} \label{eq:gaussian-est-intro}
|\wh{\mu}_N-\mu| \leq c\sigma\sqrt{\frac{\log(2/\delta)}{N}}~.
\end{equation}
Here $c$ is an absolute constant.  (In this article we focus on optimal orders of magnitude
and ignore the--important--problem of optimizing constants.)
If $X$ is indeed Gaussian, then the choice of $\wh{\mu}_N$ is simple: the empirical mean
$$
\frac{1}{N} \sum_{i=1}^N X_i
$$
has the desired accuracy at all confidence levels $\delta$.

The empirical mean also yields \eqref{eq:gaussian-est-intro} when $X$ is $L$-sub-Gaussian, that is, if for every $p \geq 2$, $\|\ol{X}\|_{L_p}\leq L \sqrt{p} \|\ol{X}\|_{L_2}$, where $\ol{X}=X-\mu$ and $\|\ol{X}\|_{L_p}=\left(\EXP|\ol{X}|^p\right)^{1/p}$,
see, for example, \cite{BoLuMa13}.

 Unfortunately, this is as far as the empirical mean takes us.
As soon as one leaves the sub-Gaussian realm, the empirical mean becomes a poor choice
and its performance deteriorates for `heavy-tailed' distributions of $X$.
In fact, for all $\delta$ there are distributions in which the estimate that follows from Chebyshev's inequality, that
\begin{equation} \label{eq:intro-1}
\PROB\left(\left|\frac{1}{N}\sum_{i=1}^N X_i - \mu \right| \geq \frac{\sigma}{ \sqrt{\delta N}} \right) \leq \delta,
\end{equation}
is sharp. In other words, while the expected value
$$
\E \left|\frac{1}{N}\sum_{i=1}^N X_i - \mu \right|
$$
is of the right order of magnitude ($\sim \sigma/\sqrt{N}$), the empirical mean exhibits rather poor concentration around $\mu$.

Thus, the empirical mean has a performance comparable to the Gaussian case only in two situations:
\begin{description}
\item{$\bullet$} For an arbitrary distribution of $X$ if one is only interested in constant confidence level (say $\delta=0.1$),
in which case the resulting accuracy is $\E|N^{-1}\sum_{i=1} X_i -\mu|$;
\item{$\bullet$} If $X$ is $L$-sub-Gaussian and one is interested in any confidence level, in which case
the error is determined by estimating the probability $\PROB(|N^{-1}\sum_{i=1}^N X_i - \mu | \geq \eta)$.
\end{description}

Perhaps surprisingly, the error one incurs in these two special and restrictive situations can be attained in full generality (though obviously the estimator one uses is not the empirical mean). One estimator that attains a ``sub-Gaussian'' performance
(i.e., an accuracy bounded by  $c N^{-1/2}\sigma \sqrt{\log(2/\delta)}$ for an absolute constant $c$) for any $X$ with finite mean and variance is the \emph{median-of-means estimator}. To compute this estimator, first the sample $X_1,\ldots,X_N$ is split into $n$ blocks $I_j$, each one of the same cardinality $m$ (here we assume without loss of generality that $n$ divides $N$). For each block $I_j$, let
$$
a_j = \frac{1}{m} \sum_{i \in I_j} X_i,
$$
and put $\wh{\mu}_N$ to be a median of $\{a_1,\ldots,a_n\}$. Setting $n \sim \log(2/\delta)$, it is straightforward to verify that this choice of $\wh{\mu}_N$ satisfies \eqref{eq:gaussian-est-intro}.
This estimator was introduced independently by Nemirovsky and Yudin \cite{NeYu83};
Jerrum, Valiant, and Vazirani \cite{JeVaVa86}; and
Alon, Matias, and Szegedy \cite{AlMaSz02}.
Another, quite different, sub-Gaussian estimator was constructed by Catoni \cite{Cat10}.

Note that unlike the empirical mean, here the procedure changes with the desired confidence.
This is indeed necessary. As it is shown by Devroye, Lerasle, Lugosi, and Oliveira \cite{DeLeLuOl16},
there is no single procedure that attains \eqref{eq:gaussian-est-intro} for all confidence levels
and for all distributions with finite second moment.

While the one-dimensional picture was well understood, in higher dimensions the situation was far less clear. Unfortunately,
establishing the 'right' notion of error in higher dimensions and with respect to a general norm can be difficult,
as parameters that are totally different in the multi-dimensional setup may `collapse' to the same object in dimension one.
However, one may still learn a lesson from the real-valued case and conclude the following:
\begin{description}
\item{$\bullet$} An estimator with accuracy of optimal order should depend on the prescribed confidence level and on the norm in question.
\item{$\bullet$} A reasonable notion of error is dictated by what happens in the two `trivial' situations---in both of which the empirical mean is essentially optimal---as in dimension one. For a real-valued random variable, when only a constant confidence is required, the error is
of the order of $\E|N^{-1}\sum_{i=1}^N X_i - \mu|$. For small values of $\delta$, the optimal accuracy is of the order of $\eta$ for which $\PROB(|N^{-1}\sum_{i=1}^N X_i - \mu | \geq \eta) \leq \delta$ when $X$ is $L$-sub-Gaussian.
The analogous objects for a random vector in $(\R^d, \| \cdot \|)$ are the expectation of the norm
\begin{equation} \label{eq:3-intro}
\E \left\|\frac{1}{N}\sum_{i=1}^N X_i - \mu\right\|
\end{equation}
and the value $\eta$ such that
\begin{equation} \label{eq:3a-intro}
\PROB \left(\left\|\frac{1}{N}\sum_{i=1}^N X_i - \mu\right\| \geq \eta\right) \leq \delta
\end{equation}
when $X$ is an $L$-sub-Gaussian random vector\footnote{Recall that $X$ is $L$-sub-Gaussian if for every $t \in \R^d$ and every $p \geq 2$, $\|\inr{X-\mu,t}\|_{L_p} \leq L \sqrt{p} \|\inr{X-\mu,t}\|_{L_2}$.}.
\end{description}
We put \eqref{eq:3-intro} in a form more convenient for us. To this end, set
$$
Y_N = \frac{1}{\sqrt{N}} \sum_{i=1}^N \eps_i (X_i - \mu)
$$
where $(\eps_i)_{i=1}^N$ are independent, symmetric, $\{-1,1\}$-valued random variables that are also independent of $(X_i)_{i=1}^N$. A standard symmetrization argument shows that
$$
\E \left\|\frac{1}{N}\sum_{i=1}^N X_i - \mu\right\|
\leq \frac{2}{\sqrt{N}} \E \left\|Y_N \right\|~.
$$
(Also observe that by the central limit theorem, $Y_N$ tends, in distribution, to the centred Gaussian random vector $G$ that has the same covariance as $X$).

As for \eqref{eq:3a-intro}, if $X$ is $L$-sub-Gaussian, then by a standard chaining argument combined with the majorizing measures theorem, one has that, with probability at least $1-\delta$,
\begin{equation} \label{eq:4-intro}
\left\|\frac{1}{N}\sum_{i=1}^N X_i - \mu\right\| \leq \frac{c(L)}{\sqrt{N}} \left(\E \|G\| + \sqrt{\log(1/\delta)} \sup_{x^* \in {\cal B}^{\circ}} \left(\E (x^*(X-\mu))^2\right)^{1/2} \right)~,
\end{equation}
where again, $G$ is the centred Gaussian vector that has the same covariance as $X$, ${\cal B}^{\circ}$ is the unit ball of the dual space\footnote{Here and in what follows we identify linear functionals on $\R^d$ with points in $\R^d$, and the action of $t \in \R^d$ is given by $x^*(x)=\inr{t,x}$, that is, the standard inner product with $t$.} to $(\R^d, \| \ \|)$, and $c(L)$ is a constant that depends on $L$ only.

Thus, if one believes that \eqref{eq:3-intro} and \eqref{eq:3a-intro} should govern the error for a general mean estimation problem in $(\R^d, \| \cdot \|)$, one arrives to the following question:

\begin{tcolorbox}

\begin{question} \label{qu:qunatitative}
Let $\| \cdot \|$, $N$ and $\delta$ be as above. Does there exist an estimator $\wh{\mu}_N$ (which may depend on $\delta$ and on the norm $\| \cdot \|$), such that, for all distributions whose covariance matrix exists, with probability at least $1-\delta$,
\begin{equation} \label{eq:main}
\|\wh{\mu}_N - \mu\| \leq \frac{c}{\sqrt{N}} \left( \max\left\{\E\|Y_N\|, \ \ \E \|G\| + R \sqrt{\log(2/\delta)} \right\}\right)~,
\end{equation}
where $c$ is an absolute constant and
$$
R=\sup_{x^* \in {\cal B}^{\circ}} \left(\E (x^*(X-\mu))^2\right)^{1/2}~?
$$
\end{question}

\end{tcolorbox}

To put Question \ref{qu:qunatitative} is some perspective, let us consider the case of the Euclidean norm  $\| \cdot \|=\| \cdot \|_2$ in $\R^d$. Let $\Tr(\Sigma)$ be the trace of the covariance matrix of $X$ and set $\lambda_1$ to be the largest eigenvalue of $\Sigma$.
Observe that
$$
\E\|Y_N\|_2 \leq (\E\|Y_N\|_2^2)^{1/2} \leq (\E \|X-\mu\|_2^2)^{1/2} = \sqrt{\Tr(\Sigma)}~,
$$
and a similar bound holds for $\E\|G\|_2$, since $Y_N$ and $G$ share the same covariance matrix. Also, because the Euclidean norm is self-dual, ${\cal B}^{\circ}=B_2^d$, the Euclidean unit ball. Therefore,
$$
R=\sup_{t \in B_2^d} \left(\E \inr{t,X-\mu}^2\right)^{1/2} \leq \sqrt{\lambda_1}~.
$$
Hence, if Question \ref{qu:qunatitative} has an affirmative answer, the resulting mean estimation error for the Euclidean norm would satisfy
\begin{equation} \label{eq:est-Euclidean}
\|\wh{\mu}_N-\mu\|_2 \leq \frac{c}{\sqrt{N}} \left(\sqrt{\Tr(\Sigma)} + \sqrt{\lambda_1 \log(2/\delta)}\right)
\end{equation}
and with probability $1-\delta$. This coincides with the performance of the empirical mean if $X$ is Gaussian (see \cite{JoLuOl17}).

As it happens, \eqref{eq:est-Euclidean} was established in \cite{LuMe16a} for an arbitrary random vector $X$ (that has a well-defined mean and covariance) using the notion of median-of-means tournaments.

In Section \ref{sec:lower} we argue that \eqref{eq:main} is not far from the best (uniform) estimate one can ever hope for.
For now simply observe that the term $N^{-1/2} R \sqrt{\log(2/\delta)}$ is truly required. Indeed, let $X$ be a Gaussian random vector with mean $\mu$. Observe that for any estimator $\wh{\psi}_N$ and any $x^* \in {\cal B}^\circ$,
$$
\|\wh{\psi}_N - \mu\| \geq |x^*(\wh{\psi}_N)-x^*(\mu)|~.
$$
Now fix $x^* \in {\cal B}^\circ$ and consider the random variable $x^*(X)$, which is a real-valued Gaussian whose mean is $x^*(\mu)$. If $\wh{\psi}_N$ performs with accuracy $\epsilon$ with probability $1-\delta$ given $X_1,\ldots,X_N$, then the real-valued estimator $x^*(\wh{\psi}_N)$ would perform with at least as good accuracy and confidence for the real-valued Gaussian variable $x^*(X)$. However, the results of \cite{Cat10} imply that the best possible accuracy for any mean estimator for a real valued Gaussian is $\sim N^{-1/2} \sigma \sqrt{\log(2/\delta)}$, and in our case, $\sigma^2 = \E (x^*(X-\mu))^2$. Taking the `worst choice' of $x^* \in {\cal B}^{\circ}$ shows that
$$
\epsilon \gtrsim \sup_{x^* \in {\cal B}^\circ} \left(\E (x^*(X-\mu))^2\right)^{1/2} \sqrt{\frac{\log(2/\delta)}{N}}=R\sqrt{\frac{\log(2/\delta)}{N}}~.
$$

Our main result is an affirmative answer to Question \ref{qu:qunatitative}, and the mean estimator that achieves
the desired accuracy is defined as follows. The estimator depends on the desired confidence $\delta\in (0,1)$ and also
on an ``accuracy parameter'' $\epsilon >0$. We show below that the procedure achieves accuracy $\epsilon$ whenever
it is at least as large as the expression on the right-hand side on \eqref{eq:main}. For simplicity of presentation we assume that
$n=\log(2/\delta)$ is an integer and that $N$ is divisible by $n$. (Otherwise an obvious modification
only effects the value of the unspecified constants so we do not lose any generality.)

\begin{tcolorbox}
\begin{description}
\item{$\bullet$} Set $\epsilon>0$.
\item{$\bullet$} Let $n = \log(2/\delta)$ and split the sample $(X_i)_{i=1}^N$ to $n$ blocks $I_j$, each of cardinality $N/n$. Set $Z_j = \frac{1}{m}\sum_{i \in I_j} X_i$.
\item{$\bullet$} Let $T$ be the set of extreme points of the dual unit ball ${\cal B}^{\circ}$. For every $x^* \in T$ set
    $$
    S_{x^*} = \left\{ y \in \R^d : |x^*(Z_j) - x^*(y)| \leq \epsilon \right\} \ {\rm for \ more \ than \ } \frac{n}{2} \ {\rm blocks}.
    $$
    \item{$\bullet$} Set $\mathbbm{S}(\epsilon) = \bigcap_{x^* \in T} S_{x^*}$ and select $\wh{\mu}_N(\epsilon,\delta)$ to be any point in $\mathbbm{S}(\epsilon)$.
\end{description}
\end{tcolorbox}

Note that $S_{x^*}$ is a union of intersections of shifts of the same `slab' in $\R^d$, defined by the linear functional $x^*$ and of `width' $\epsilon$. Thus, each intersection is just a (data dependent) slab, making $S_{x^*}$ to be the union of slabs defined by $x^*$. As a result, $\mathbbm{S}(\epsilon)$ is an intersection of unions of slabs generated by the extreme points of the dual unit ball of the given norm.

Our main result is the following---formulated using the notation introduced previously.
\begin{theorem} \label{thm:main}
There exist absolute constants $c,c^\prime$ such that the following holds.
Given a norm $\| \ \|$, confidence parameter $\delta\in (0,1)$ and sample size $N$, if
\begin{equation} \label{eq:main1}
\epsilon \geq \frac{c}{\sqrt{N}} \left( \max\left\{\E\|Y_N\|, \ \ \E \|G\| + R \sqrt{\log(2/\delta)} \right\}\right)~,
\end{equation}
then the estimator $\wh{\mu}_N(\epsilon,\delta)$ defined above satisfies that, with probability at least $1-c^\prime\delta$, $\mathbbm{S}(\epsilon)$ is nonempty, and
$$
\|\wh{\mu}_N(\epsilon,\delta) - \mu\| \leq \epsilon~.
$$
\end{theorem}

\begin{remark}
Observe that Theorem \ref{thm:main} also implies that if $\epsilon$ is as in \eqref{eq:main1} then the set $\mathbbm{S}(\epsilon)$ is bounded. Moreover, the set is also closed because each $S_{x^*}$ is closed.
\end{remark}

The estimator $\wh{\mu}_N(\epsilon,\delta)$ has the disadvantage that it requires the knowledge of the accuracy level $\epsilon$.
However, the achievable optimal accuracy depends on the distribution and it is generally unknown to the statistician.
Luckily, it is easy to use the theorem above to construct an estimator that does not depend on such previous knowledge and yet
achieves the same performance bound. We may simply define our estimator $\wt\mu_N=\wt\mu_N(\delta)$ as follows.
Let $\epsilon_0=\inf\{\epsilon>0: \mathbbm{S}(\epsilon)\neq \emptyset \}$. The sets $\mathbbm{S}(\epsilon)$ for $\epsilon>\epsilon_0$
are nested and compact and therefore $\cap_{\epsilon>\epsilon_0}\mathbbm{S}(\epsilon)\neq\emptyset $. We define  $\wt\mu_N$ to be an arbitrary element of $\cap_{\epsilon>\epsilon_0}\mathbbm{S}(\epsilon)$.
It follows from Theorem \ref{thm:main} that for $\epsilon$ satisfying \eqref{eq:main1}, with probability at least $1-\delta$ , $\mathbbm{S}(\epsilon)\neq \emptyset$, and in particular, $\wt\mu_N\in \mathbbm{S}(\epsilon)$.  Hence, we obtain the following.
\begin{corollary}
There exist absolute constants $c,c^\prime$ such that the following holds.
Given a norm, confidence parameter $\delta\in (0,1)$ and sample size $N$, if
\begin{equation} \label{eq:main2}
\epsilon \geq \frac{c}{\sqrt{N}} \left( \max\left\{\E\|Y_N\|, \ \ \E \|G\| + R \sqrt{\log(2/\delta)} \right\}\right)~,
\end{equation}
then the estimator $\wt{\mu}_N$ satisfies that, with probability at least $1-c^\prime \delta$,
$$
\|\wt{\mu}_N - \mu\| \leq \epsilon~.
$$
\end{corollary}

Theorem \ref{thm:main} is established using a general fact that is of independent interest: we construct an effective \emph{uniform} median-of-means estimator in a class of real valued functions, as described in the next section.

\subsection*{Related work}

The multivariate median-of-means estimators that behave well under heavy-tailed distributions
have been the subject of intensive study.
Minsker \cite{Min15} and Hsu and Sabato \cite{HsSa16} defined and analyzed multivariate extensions of the
median-of-means estimator,
see also Lerasle and Oliveira \cite{LerasleOliveira_Robust}.
The first truly sub-Gaussian estimator (under the Euclidean norm) was shown to exist by Lugosi and Mendelson \cite{LuMe16a}.
See Joly, Lugosi, and Oliveira \cite{JoLuOl17} for an earlier attempt and
Catoni and Giulini \cite{CaGi17} for a different estimator.

Minsker \cite{Min16} and Catoni and Giulini \cite{CaGi17} consider estimating the mean of random matrices
based on an i.i.d.\ sample under the spectral norm and the Hilbert-Schmidt norm. They both prove sub-Gaussian
performance bounds but the bounds of these papers fall short, in various aspects, of the optimal order of magnitude achieved by the
estimator of Theorem \ref{thm:main} above. As far as we know, estimators achieving the accuracy/tradeoff of Theorem \ref{thm:main} have only been known for the Euclidean norm.

\section{Uniform median-of-means estimators}

In this section we explore the next problem:

\begin{tcolorbox}
Let $F$ be a class of functions on a probability space $(\Omega,\nu)$ and let $\delta \in (0,1)$. Given an independent sample $(X_1,\ldots,X_N)$ distributed according to $\nu^N$, find an estimator $\wh{\Phi}_N$, such that, with probability at least $1-\delta$, for every $f \in F$, $|\wh{\Phi}_N(f) - \E f |$ is small.
\end{tcolorbox}

The obvious choice of $\wh{\Phi}_N$ is simply the standard median-of-means estimator we use for a single random variable. However, expecting $\wh{\Phi}_N$ to have the `individual' sub-Gaussian error is too optimistic.
The best uniformly achievable accuracy must depend on some appropriate notion of the `size' of the class $F$.

To address the problem above, fix integers $n$ and $m$ and let $N=mn$. As before, we split the given sample to $n$ blocks, each one of cardinality $n$, while keeping in mind that the natural choice is $n \sim \log(2/\delta)$.
Our goal is to find the smallest possible value of $r$ such that $\sup_{f\in F} |\wh{\Phi}_N(f) - \E f |\le r$ with probability at least
$1-\delta$.

Recall that if one would like to ensure that the median-of-means estimator performs with an error of at most $r$ for a \emph{single function} $f \in F$, then it suffices that
\begin{equation} \label{eq:single}
\PROB\left(\left|\frac{1}{m}\sum_{i=1}^m f(X_i)-\E f \right|\geq r \right) \leq \frac{1}{2}-\theta
\end{equation}
for some $\theta>0$. Indeed, if \eqref{eq:single} holds then with probability at least $1-2\exp(-c\theta^2n)$,
$$
\E f  - r \leq \frac{1}{m}\sum_{i \in I_j} f(X_i) \leq \E f+r
$$
for more than $n/2$ of the blocks $I_j$, where $c$ is an absolute constant. However, a uniform result calls for a little more flexibility. Firstly, there is a need to have a larger number of `good'  blocks $I_j$. It suffices that for any fixed function one controls $0.9n$ of them. Clearly, that may be achieved if \eqref{eq:single} holds for $1/2 - \theta \leq 0.05$.  With that in mind, let
$$
p_m(\eta) = \sup_{f \in F} \PROB\left(\left|\frac{1}{m}\sum_{i=1}^m f(X_i)-\E f \right|\geq \eta \right)~.
$$
From here on we write at times $p_m$ instead of $p_m(\eta)$. We set $D$ to be the unit ball in $L_2(\nu)$ and let ${\cal M}(F,rD)$ be the maximal cardinality of a subset of $F$ that is $r$-separated with respect to the $L_2(\nu)$ norm. We also denote $F-F=\{f_1-f_2 : f_1,f_2 \in F\}$.

Let us describe the performance of the uniform median-of-means estimator:
\begin{theorem} \label{thm:uniform-MOM}
There exist absolute constants $c_0,\ldots,c_4$ for which the following holds. Set $\eta_0,\eta_1$ and $\eta_2 \geq c_0 \eta_1/\sqrt{m}$ that satisfy the following:
\begin{description}
\item{$(1)$} $p_m(\eta_0) \leq 0.05$~;
\item{$(2)$} $\log {\cal M}(F,\eta_1 D) \leq c_2 n \log(e/p_m(\eta_0))$~;
\item{$(3)$} $\E \sup_{w \in \ol{W}} \left|\sum_{i=1}^N \eps_i w(X_i) \right| \leq c_3 \eta_2 N$~,
\end{description}
where $W=(F-F) \cap \eta_1 D$ and $\ol{W}=\{w -\E w : w \in W\}$.

Let $r=\eta_0+\eta_2$. Then with probability at least $1-2\exp(-c_4n)$, for any $f \in F$ one has that
$$
\left| \frac{1}{m} \sum_{i \in I_j} f(X_i) - \E f  \right| \leq r \ \ {\rm for \ at \ least \ } 0.6n \ {\rm blocks \ } I_j~.
$$
\end{theorem}
To put Theorem \ref{thm:uniform-MOM} in some perspective, note that
$\eta_0$ captures the worst individual error caused by a function in $F$. Moreover, as noted previously, the standard median-of-means estimator would perform with accuracy $\eta_0$ and confidence $1-\delta$ if
$$
\PROB\left(\left|\frac{1}{m}\sum_{i=1}^m f(X_i)-\E f \right|\geq \eta_0 \right) \leq 0.05,
$$
and by Chebyshev's inequality, one may set
$$
\eta_0 \gtrsim \left(\E(f(X)-\E f)^2\right)^{1/2} \cdot \frac{1}{\sqrt{m}} \sim \left(\E(f(X)-\E f)^2\right)^{1/2} \cdot \sqrt{\frac{\log(2/\delta)}{N}},
$$
as one would expect from a sub-Gaussian estimate.

In contrast, the role of $\eta_2$ is to calibrate the impact of the `size' of $F$.


\proof Fix $f \in F$ and let $\delta_j$ be the indicator of the event $\left|\frac{1}{m}\sum_{i \in I_j} f(X_i)-\E f \right|\geq \eta_0$.
By a standard binomial tail estimate, for $k \geq 0.06n$,
$$
\PROB\left( \left|\{j : \delta_j = 1\}\right| \leq k \right) \geq 1-2\exp(-c_0 k \log(ek/p_m n))~.
$$
In particular, $|\{j : \delta_j =1\}| \leq 0.1n$ with probability at least $1-2\exp(-c_1n \log(e/p_m))$.

The importance of the high-probability estimate is seen in the next step of the proof: one may control all the elements of an $\eta_1$-net of $F$ (with respect to the $L_2(\nu)$ norm) as long as its cardinality is at most $\exp(c_2 n \log(e/p_m))$. Indeed, by the union bound, with probability at least $1-2\exp(-c_3 n \log(e/p_m))$, for every $h$ in the net there are at least $0.9n$  blocks $I_j$ such that
$$
\left|\frac{1}{m} \sum_{i \in I_j} h(X_i) - \E h \right| \leq \eta_0~.
$$
The final and crucial step in the proof is passing from the net to the entire class: for every $f \in F$ set $\pi f$ to be the best approximation to $f$ in the net. Thus, $\|f-\pi f\|_{L_2} \leq \eta_1$. We show that for every $f \in F$ there are at most $0.2 n$  blocks $I_j$ such that
\begin{equation} \label{eq:osc-in-proof}
\left|\frac{1}{m}\sum_{i \in I_j} (f(X_i)-\E f) -  \frac{1}{m}\sum_{i \in I_j} (\pi f - \E \pi f)(X_i)\right| \leq \eta_2.
\end{equation}
If that is indeed the case then for every $f \in F$ there are at least $0.7n$  blocks for which
\begin{eqnarray*}
\lefteqn{
 \left|\frac{1}{m} \sum_{i \in I_j} f(X_i) - \E f\right|   }
\\
& \leq & \left|\frac{1}{m} \sum_{i \in I_j} (\pi f)(X_i) - \E \pi f \right| + \left|\frac{1}{m}\sum_{i \in I_j} (f-\E f)(X_i) - \frac{1}{m}\sum_{i \in I_j} (\pi f - \E \pi f) (X_i)\right|
\\
& \leq & \eta_0 + \eta_2~,
\end{eqnarray*}
as required.

It remains to prove \eqref{eq:osc-in-proof}. To this end, note that $f - \pi f \in (F-F) \cap \eta_1 D = W$, and thus
$f - \E f - (\pi f - \E \pi f ) \in \ol{W}$ were $\ol{W}=\{w - \E w : w \in W\}$.
Hence, the proof is completed once it is established that, with probability at least $1-2e^{-c_4n}$,
\[
S \defeq \sup_{w \in \ol{W}} \left| \left\{ j : \left|\frac{1}{m}\sum_{i \in I_j} w(X_i) \right| \geq \eta_2 \right\} \right| \leq 0.2n~.
\]
To control $S$, note that by the bounded differences inequality (see, e.g., \cite{BoLuMa13}) there is an absolute constant $c_1$ such that
$$
\PROB( S \geq \E S+ 0.1n ) \leq 2\exp(-c_1n)~.
$$
Thus, all that remains is to show that $\E S \leq 0.1n$. Observe that for any $(a_j)_{j=1}^n$,
$$
|\{j : |a_j| \geq \eta\}| = \sum_{j=1}^n \IND_{\{|a_j| \geq \eta\}} \leq \frac{1}{\eta} \sum_{j=1}^n |a_j|~.
$$
Hence, by standard methods of empirical processes,
via an analogous argument to that in \cite{LuMe16a}, one has
\begin{eqnarray*}
\lefteqn{
\E \sup_{w \in \ol{W}} \left|\left\{j : \left|\frac{1}{m}\sum_{i \in I_j} w(X_i)\right| \geq \eta_2 \right\} \right|
} \\
& \leq & \frac{1}{\eta_2} \E \sup_{w \in \ol{W}} \sum_{j=1}^n \left|\frac{1}{m}\sum_{i \in I_j} w(X_i)\right|
\\
& \leq & \frac{1}{\eta_2} \E \sup_{w \in \ol{W}} \sum_{j=1}^n \left( \left|\frac{1}{m}\sum_{i \in I_j} w(X_i)\right| - \E \left|\frac{1}{m}\sum_{i \in I_j} w(X_i)\right| \right) + \frac{n}{\eta_2}\sup_{w \in \ol{W}} \E \left|\frac{1}{m}\sum_{i \in I_j} w(X_i)\right|
\\
& \leq & \frac{2}{\eta_2} \E \sup_{w \in \ol{W}} \left| \sum_{j=1}^n \eps_j \left(\frac{1}{m} \sum_{i \in I_j} w(X_i) \right) \right| + \frac{n}{\eta_2} \cdot \sup_{w \in \ol{W}} \frac{\|w\|_{L_2}}{\sqrt{m}}
\\
& \leq & \frac{4n}{\eta_2} \left(\E \sup_{w \in \ol{W}} \left|\frac{1}{N} \sum_{i=1}^N \eps_i w(X_i) \right| + \frac{\eta_1}{\sqrt{m}}\right)~.
\end{eqnarray*}
In particular, $\E S \leq 0.1n$ provided that
\[
\E \sup_{w \in \ol{W}} \left| \sum_{i=1}^N \eps_i w(X_i)\right| \leq c_2 \eta_2 N \ \ \ \ {\rm and} \ \ \ \ \frac{\eta_1}{\sqrt{m}} \leq c_3 \eta_2~,
\]
as we assumed.
\endproof
\begin{remark}
Note that if $F$ is a finite class and $\log |F| \leq c_2 n \log(e/p_m(\eta_0))$ then $\wh{\Phi}_N$ performs with accuracy $\eta_0$. The proof follows from the standard bound on the performance of the median-of-means estimator for each real random variable $f(X)$ and a straightforward application of the union bound.
\end{remark}

\section{Estimation with respect to a general norm}

In this section we establish Theorem \ref{thm:main} by invoking Theorem \ref{thm:uniform-MOM}.

Let $\| \cdot \|$ be a norm on $\R^d$ and let ${\cal B}^\circ$ be the unit ball of the dual norm. Recall that for any $v \in \R^d$,
$$
\|v\| = \sup_{x^* \in {\rm ext}({\cal B}^{\circ})} x^*(v)~,
$$
where ${\rm ext}({\cal B}^\circ)$ denotes the set of extreme points in ${\cal B}^\circ$, and that the empirical average within block $I_j$, for $1 \leq j \leq n$, is denoted by
$$
Z_j = \frac{1}{m} \sum_{i \in I_j} X_i~.
$$
Let $r>0$ be as in Theorem \ref{thm:uniform-MOM} for the class of functions $F = \left\{ x^*(\cdot) : x^* \in {\rm ext}({\cal B}^\circ)  \right\}$ and with the respect to the measure $\nu$ endowed by $X-\mu$.
Finally, let ${\cal A}$ be the event for which the assertion of Theorem \ref{thm:uniform-MOM} holds.

Consider the sets 
$$
S_{x^*}=\left\{ y \in \R^d : \left|x^*(Z_j)-x^*(y)\right| \leq r \ \ {\rm for \ more \ than \ } \frac{n}{2} \ {\rm indices} \ j \right\}
$$
and put $\wh{\mu}_N(\epsilon,\delta)$ to be any point that belongs to the set
\begin{equation} \label{eq:slabs-intersection}
\mathbbm{S}(\epsilon)=\bigcap_{x^* \in {\rm ext}({\cal B}^\circ)} S_{x^*}~.
\end{equation}

To show that selecting $\wh{\mu}_N(\epsilon,\delta) \in \mathbbm{S}(\epsilon)$ has the desired properties, fix a sample $(X_i)_{i=1}^N \in {\cal A}$.
First, observe that $\mathbbm{S}(\epsilon)$ is nonempty as it contains  $\mu$. Indeed,  setting $f(x)=x^*(x)$, it is evident that
$$
\E f(X-\mu) = 0 \ \ \ \ {\rm  and} \ \ \ \ \frac{1}{m}\sum_{i \in I_j} f(X_i-\mu) = x^*(Z_j)-x^*(\mu)~.
$$
By Theorem \ref{thm:uniform-MOM} it follows that
$$
| x^*(Z_j) - x^*(\mu)| \leq r
$$
for a majority of the indices $j$, which means that $\mu \in S_{x^*}$ for every $x^* \in {\rm ext}({\cal B}^\circ)$.

Next, one has to show that if $y \in \mathbbm{S}(\epsilon)$, then $\|y - \mu \|$ is `small'. To that end, observe that for every $x^* \in {\rm ext}({\cal B}^\circ)$ there is some index $j$ such that
$$
\left|x^*(Z_j)-x^*(y)\right| \leq r \ \ \ {\rm and} \ \ \ \left|x^*(Z_j)-x^*(\mu)\right| \leq r~,
$$
because both conditions hold for more than half of the indices $j$. Thus,
$$
\left|x^*(y)-x^*(\mu)\right| \leq \left|x^*(Z_j)-x^*(y)\right|+\left|x^*(Z_j)-x^*(\mu)\right| \leq 2r~.
$$
Finally, recalling that $\|v\|=\sup_{x^* \in {\rm ext}({\cal B}^\circ)} x^*(v)$, one has that
$$
\|y-\mu\| = \sup_{x^* \in {\rm ext}({\cal B}^\circ)} |x^*(y)-x^*(\mu)| \leq 2r~,
$$
as claimed.

To complete the proof of Theorem \ref{thm:main} let us bound $\eta_0$ and $\eta_2$. To that end, recall that
$$
R = \sup_{x^* \in {\cal B}^\circ} \left(\E (x^*(X-\mu))^2\right)^{1/2}~,
$$
$G$ is the centred Gaussian vector that has the same covariance as $X$, and
$$
Y_N = \frac{1}{\sqrt{N}} \sum_{i=1}^N \eps_i (X_i -\mu)~.
$$
%


We show that the three conditions of Theorem \ref{thm:uniform-MOM} can be controlled when $F=\{x^*(\cdot) : x^* \in {\rm ext}({\cal B}^\circ)\}$ and with respect to the measure $\nu$ endowed by $\ol{X}=X-\mu$.

To verify $(1)$, fix $x^* \in {\cal B}^{\circ}$ and note that
$$
\PROB\left( \left|\frac{1}{m} \sum_{i=1}^m x^*(X_i-\mu) \right| \geq \eta_0\right) \leq \frac{\E|x^*(X-\mu)|^2}{\eta_0^2 m} = \frac{\E|x^*(X-\mu)|^2 \log(2/\delta)}{\eta_0^2 N}
$$
where we have used the fact that $n = \log(2/\delta)$.
(Recall that we assume, without loss of generality, that $\log(2/\delta)$ is an integer that divides $N$.)
Thus, to ensure that $p_m \leq 0.05$  it suffices that
\[
\eta_0 \geq c_0 R \sqrt{\frac{\log(2/\delta)}{N}}~.
\]
Turning to $(2)$, we identify ${\cal B}^\circ$ with the set $\{t \in \R^d : \sup_{x \in {\cal B}} \inr{t,x} \leq 1\}$, and the action of a functional $x^*$ associated with $t$ is given by $x*(x)=\inr{t,x}$. We also abuse notation and denote by $D$ the unit ball of the $L_2(\ol{X})$ norm endowed on $\R^d$ by identifying each $t \in \R^d$ with a linear functional. By Sudakov's inequality (see \cite{LeTa91}), there is an absolute constant $c$ such that
$$
\log {\cal M}({\cal B}^\circ, \eta_1 D) \leq c\eta_1^{-2}(\E \sup_{x^* \in {\cal B}^\circ} x^*(G))^2= c \left(\frac{\E \|G\|}{\eta_1}\right)^2~,
$$
implying that one may set
\[
\eta_1 =c_1 \frac{\E \|G\|}{\sqrt{n}}~.
\]
In particular, this forces the constraint
\[
\eta_2 \geq c_2 \frac{\E \|G\|}{\sqrt{N}}~.
\]

Finally, to control $(3)$, observe that
$$
W \subset \{\inr{t,\cdot} : t \in 2{\cal B}^\circ \cap \eta_1 D\}~.
$$

Therefore, one has to show that
$$
\E \sup_{t \in 2{\cal B}^\circ \cap \eta_1 D} \left|\frac{1}{\sqrt{N}}\sum_{i=1}^N \eps_i \inr{t,(X_i-\mu)} \right| = \E\sup_{t \in 2{\cal B}^\circ \cap \eta_1 D} |\inr{t,Y_N}| \le \eta_2 \sqrt{N}~,
$$
where $Y_N=N^{-1/2}\sum_{i=1}^N \eps_i (X_i-\mu)$. Clearly, it suffices that
\begin{equation} \label{eq:cond-eta-2-general-1}
\eta_2 \geq c_3 \frac{\E \|Y_N\|}{\sqrt{N}}~,
\end{equation}
and one may set
$$
\eta_2=\frac{c_4}{\sqrt{N}} \max\left\{\E \|Y_N\|, \E \|G\| \right\}~.
$$
Now Theorem \ref{thm:main} follows from Theorem \ref{thm:uniform-MOM}.
\endproof

\begin{remark}
Note that the choices of $\eta_0,\eta_1$ and $\eta_2$ need not be optimal for each $F$ and $\ol{X}$ as above. Indeed, $\eta_1$ was chosen via Sudakov's inequality which is not always sharp, and $\eta_2$ was determined after the `localization' $2{\cal B}^\circ \cap \eta_1 D$ was replaced by $2{\cal B}^\circ$. Therefore, it stands to reason that there are cases in which the resulting estimate may be improved with more care. However, as we explain in the next section, Theorem \ref{thm:main} is likely to be the best uniform result that one can hope for.
\end{remark}

\section{Lower bounds} \label{sec:lower}

In this section we discuss the optimality of the upper bound of Theorem \ref{thm:main}.
In the introduction we already pointed out that  the term $N^{-1/2} R \sqrt{\log(2/\delta)}$ is inevitable.
Here we discuss the necessity of the term $N^{-1/2} \E \|G\|$ (and, equivalently, the term $N^{-1/2} \E\|Y_N\|$ since
$\lim_{N\to \infty} \E\|Y_N\|= \EXP\|G\|$ by the central limit theorem).

Here we show that the order of
magnitude of the bound of Theorem \ref{thm:main} is essentially un-improvable even if one only considers
isotropic Gaussian distributions, with one minor caveat.
Recall that the proof of Theorem \ref{thm:main} uses Sudakov's inequality to ensure that
\begin{equation} \label{eq:entropy-optimal}
\log {\cal M}({\cal B}^\circ, \eta_1 D) \lesssim n~,
\end{equation}
and then the contribution to the error is $\sim \eta_1/\sqrt{m}=\eta_1 \sqrt{n/N}$.
As noted previously, while it is convenient to use Sudakov's inequality, its application may be loose.
A more accurate upper estimate on the error is $\sim \eta_1/\sqrt{m}$ where $\eta_1$
is the smallest value for which \eqref{eq:entropy-optimal} holds.

We show now that if $X$ is a Gaussian measure whose covariance is the identity matrix,
then this more accurate upper estimate is actually a lower bound as well.

To formulate the lower bound, let $G$ be the standard Gaussian random vector in $\R^d$ and denote by $\nu$ the corresponding measure. We study the performance of an arbitrary mean estimation procedure with respect to a norm $\| \cdot \|$, for a collection of Gaussian measures endowed by $\{G + t : t \in T\}$ and where $T \subset \R^d$ is a well chosen set. Note that for any $X = G+t$ one has that $\ol{X}=G$, let $D \subset \R^d$ be the unit ball endowed by the norm $L_2(\ol{X})=L_2(\nu)$ (which in this case is simply $B_2^d$), and set ${\cal B}$ to be the unit ball of $(\R^d, \| \cdot \|)$.

\begin{theorem} \label{thm:lower}
There exist absolute constants $c_1$ and $c_2$ for which the following holds. Let $n \leq N$ and assume that $\log {\cal M}({\cal B}^\circ, \eta D) \geq c_1 n$. There is a set $T \subset \R^d$ such that any mean estimator $\wh{\Psi}_N$ that performs with confidence $1/2$ with respect to all the Gaussian measures $\{G + t : t \in T\}$, cannot perform with higher accuracy than $c_2\eta \sqrt{n/N}$.
\end{theorem}

\begin{remark}
An immediate outcome of Theorem \ref{thm:lower} is that when Sudakov's inequality is sharp at scale $\eta$,
the lower bound on the accuracy that holds with constant confidence is indeed $\sim \E\|G\|/\sqrt{N}$.
\end{remark}

\proof
To define the set $T$, observe that if $\log {\cal M}({\cal B}^\circ, \eta D) \geq c_1 n$ then by the duality theorem of metric entropy \cite{ArMiSz04}, and since $D^\circ = B_2^d$, it follows that $\log {\cal M}(D, c_2\eta {\cal B}) \geq c_3 n$. In other words, the set $D$ contains a subset that is $c_2\eta$-separated with respect to the norm $\| \ \|$ and whose cardinality is $c_3n$. 
Set $R=\sqrt{n/N}=1/\sqrt{m}$ and $r=c_2\eta/\sqrt{m}$. Clearly,
$$
\log {\cal M}(D, c_2\eta {\cal B}) = \log {\cal M}(D, (r/R) {\cal B}) = \log {\cal M}(R D, r {\cal B}) \geq c_3n~,
$$
and let $T \subset RD$ be the $r$-separated set with respect to the norm $\| \ \|$.

Now, assume that there is a mean estimator that performs with confidence $1/2$ and accuracy $r/3$ for every one of the Gaussian random vector $G+t$, $t \in T$, and let as reach a contradiction.

If we denote by $\wh{\Psi}_N$ such an estimator, it follows that for every $t \in T$ there is a set ${\cal A}_t \subset (\R^d)^N$ with $\nu^N(A_t) \geq 1/2$, such that for $\omega \in {\cal A}_t + (t,\ldots,t)$ one has that $\|\wh{\Psi}_N(\omega) - t \| \leq r/3$. Moreover, since the set $T$ is $r$-separated, it is evident that the sets $U_t={\cal A}_t + (t,\ldots,t)$ are disjoint. Indeed, if $\omega \in ({\cal A}_x + (x,\ldots,t)) \cap ({\cal A}_y + (y,\ldots,y))$ then $\wh{\Psi}_N$ would have to be `close' to both $x$ and $y$, which is impossible.


Observe that $\nu^N$ is the standard Gaussian measure on $(\R^d)^N$, and $\|(t,\ldots,t)\|_{L_2(\nu^N)} =\sqrt{N}\|t||_2 \leq \sqrt{N}R$.
Using the same argument
as in Talagrand's proof of the dual Sudakov inequality \cite{LeTa91} (see also \cite{MR3718154}) and recalling that $t \in RD=RB_2^d$,
$$
\nu^N(U_t) \geq \nu^N({\cal A}_t)\exp(-cN\|t\|_2^2) \geq \frac{1}{2}\exp(-cNR^2)~,
$$
for an absolute constant $c$.

On the other hand, the sets $U_t$ are disjoint and thus
$$
1 \geq \nu^N\left(\bigcup_{t \in T} U_t\right) = \sum_{t \in T} \nu^N(U_t) \geq c|T| \exp(-c^\prime NR^2)~.
$$
In particular,
$$
\log |T| \leq c^{\prime \prime} NR^2 = c^{\prime \prime} n~,
$$
which is a contradiction to our choice of $T$.
\endproof

\bibliographystyle{plain}

\end{document}